IMS Collections
**Pushing the Limits of Contemporary Statistics: Contributions in Honor of Jayanta K. Ghosh**
Vol. 3 (2008) 60–70

# Data-dependent probability matching priors for empirical and related likelihoods

**Rahul Mukerjee**[1]

*Indian Institute of Management Calcutta*

**Abstract:** We consider a general class of empirical-type likelihoods and develop higher order asymptotics with a view to characterizing members thereof that allow the existence of possibly data-dependent probability matching priors ensuring approximate frequentist validity of posterior quantiles. In particular, for the usual empirical likelihood, positive results are obtained. This is in contrast with what happens if only data-free priors are entertained.

## Contents

## 1. Introduction

Although empirical and related likelihoods have received significant attention (see [11, 14] and the references therein), their study from a Bayesian perspective began only in recent years. Lazar [10] pioneered an investigation on the validity of the empirical likelihood for posterior inference and examined, mostly by simulation, the frequentist properties of posterior empirical likelihood intervals. In another significant development, Schennach [18] proposed a Bayesian exponentially tilted empirical likelihood arising as a nonparametric limit of a Bayesian procedure which places a kind of noninformative prior on the space of distributions. Starting from a general class of empirical-type likelihoods for the population mean, Fang and Mukerjee [7] characterized its members which admit probability matching priors in the sense of allowing posterior credible sets with approximate frequentist validity.

Along the line of what is traditionally done in parametric inference based on the true likelihood, Fang and Mukerjee [7] entertained only priors that are free from the data. They observed, among other things, that none of the standard empirical-type likelihoods that have been proposed and widely studied in the literature, including the usual empirical likelihood, admits a probability matching prior even with margin of error $o(n^{-1/2})$, where $n$ is the sample size. This is somewhat disappointing and, given the popularity of these likelihoods, prompts one to investigate the consequences of working with possibly data-dependent priors with the hope that this may yield more positive results. The present article aims at exploring this issue with reference to the same general class as in [7]. Satisfyingly, it is seen that at least for the usual empirical likelihood, positive results then emerge even with margin of error $o(n^{-1})$.

[1]Joka, Diamond Harbour Road, Kolkata 700 104, India, e-mail: rmuk1@hotmail.com

Probability matching priors have been studied extensively in parametric inference – see, for example, [6, 9, 12], and [20]. A key difference with the parametric case is that here we are not working with the true density-based likelihood and, as such, a shrinkage argument, suggested originally by J. K. Ghosh to the present author, that simplifies the frequentist calculations there ([6], Ch. 1) is no longer applicable. As a result, one has to employ a direct Edgeworth expansion based on approximate cumulants.

While the present work seems to be the first attempt towards exploring the higher-order asymptotics on data-dependent priors with empirical-type likelihoods, such priors have received considerable attention in recent years in parametric inference via the true density-based likelihood. A brief indication of this literature may interest the reader. A key reference in this regard is [21], where it was found that for certain mixture models, no data-free improper prior yields a proper posterior and no data-free proper prior entails frequentist validity of posterior quantiles with margin of error $O(n^{-1})$, while both problems are solved by a data-dependent prior. Furthermore, such data-dependent priors were shown to approximate data-free priors, in addition to enjoying desirable properties like asymptotic minimaxity. Prior to [21], data-dependent priors were considered, among others, by [15] and [17] in the context of mixture models with an unknown number of components. Sweeting [19] investigated the crucial role played by data-dependent probability matching priors when the sample size is stochastic, as happens, for instance, with censoring or a stopping rule. Reid et al. [16] reviewed and discussed a notion of strong matching which requires data-dependent priors.

In the context of parametric inference, Clarke and Yuan [4] studied partial information reference priors obtained through the maximization of conditional Shannon mutual information. These priors are often data-dependent in the sense of involving statistics that are associated with nuisance parameters and capture helpful side information. The information theoretic interpretation for these priors was also discussed at length by [4]. Clarke ([3], Subsection 7.2) discussed data-dependent priors in the light of the Freedman–Purves Theorem [8], which often forms the basis of the argument put forward by orthodox Bayesians against such priors on the grounds of incoherence. He argued that the implications of this theorem are narrower than commonly appreciated, suggested a remedy in the form of a criterion of information boundedness, and observed that the data-dependent priors in [4], [16], and [21] are, indeed, information bounded.

The interested reader may refer to the papers cited in the last two paragraphs for further references on data-dependent priors in parametric inference.

## 2. A general class of empirical-type likelihoods: posterior quantiles

Let $X_1, \ldots, X_n$ be independent scalar-valued random variables from an unknown common distribution with an unknown mean $\theta$. The parameter space for $\theta$ is the real line or an open interval thereof. The $X_i$ are supposed to be absolutely continuous and the first four population moments are assumed to exist [2]. These assumptions justify an Edgeworth expansion used later. Let $\bar{X} = n^{-1}\sum_{i=1}^{n} X_i$, $m_s = n^{-1}\sum_{i=1}^{n}(X_i - \bar{X})^s$, $s = 2, 3, \ldots$, $g_3 = m_3/m_2^{3/2}$, $g_4 = m_4/m_2^2$ and $y = y(\theta) = (n/m_2)^{1/2}(\theta - \bar{X})$. Write $\phi(\cdot)$ for the standard univariate normal density.

As in [7], we consider a general class of empirical-type likelihoods of the form

$$(1)\qquad \begin{aligned} L(\theta) \quad \propto \quad & \phi(y)[1 + n^{-1/2}\{a_1(g_3)y + a_3(g_3)y^3\} \\ & +n^{-1}\{b_0(g_3, g_4) + b_2(g_3, g_4)y^2 + b_4(g_3, g_4)y^4 + b_6(g_3, g_4)y^6\} \\ & +o_p(n^{-1})], \end{aligned}$$

where the $a_i(\cdot)$ are polynomials in $g_3$ and the $b_i(\cdot)$ are polynomials in $g_3$ and $g_4$, the coefficients therein being constants free from $n$. These polynomials depend on the particular likelihood. Note that $L(\theta)$ depends on $\theta$ through $y = y(\theta)$ but does not involve any other population parameter, and that $y$ is the standard pivotal quantity for inference on $\theta$ when the population variance is unknown. Furthermore, the terms of order $n^{-1/2}$ and $n^{-1}$ in (1) aim at taking care of the unknown skewness and kurtosis of the population via their sample analogs $g_3$ and $g_4$ respectively. Indeed, the class (1) is very wide and, as discussed later in Section 4, covers all major empirical-type likelihoods proposed in the literature.

With reference to any likelihood in the class (1), we consider a possibly data-dependent prior of the form

$$(2)\qquad \pi(\theta) = \exp\{\psi(\theta, m_2, g_3)\},$$

where $\psi(\cdot)$ is a smooth function with functional form free from $n$. We aim at characterizing the $a_i(\cdot)$ and $b_i(\cdot)$ in (1) so as to allow the existence of a prior of the form (2) that entails frequentist validity of the posterior quantiles of $\theta$, with margin of error $o(n^{-1/2})$ or $o(n^{-1})$. The class (2) can be motivated as follows. For any empirical-type likelihood as in (1) and a data-free prior, up to the first order of approximation, the quantity $(\theta - \bar{X})m_2^{-1/2}$ represents a standardized version of $\theta$ in the posterior setup; see [7]. This prompts one to consider a data-dependent prior of the form

$$(3)\qquad \pi(\theta) = \exp\{(\theta - \bar{X})m_2^{-1/2}\chi(g_3)\},$$

where the multiplier $\chi(g_3)$ is a smooth function of $g_3$, with functional form free of $n$, that aims at taking care of the population skewness (an attempt to take care of the population kurtosis would involve a more elaborate data dependent prior and a discussion of this is deferred till Section 5). Clearly, the prior in (3) is equivalent to

$$(4)\qquad \pi(\theta) = \exp\{\theta m_2^{-1/2}\chi(g_3)\},$$

because they both lead to the same posterior. The exponent in (4) is a smooth function of $\theta$, $m_2$ and $g_3$, and the class (2) incorporates all priors that share this feature.

Let $\psi_i(t_1,t_2,t_3)=\partial\psi(t_1,t_2,t_3)/\partial t_i$, $\psi_{ij}(t_1,t_2,t_3)=\partial^2\psi(t_1,t_2,t_3)/\partial t_i\partial t_j$, $\psi_i=\psi_i(\bar{X},m_2,g_3)$ and $\psi_{ij}=\psi_{ij}(\bar{X},m_2,g_3)$, $i,j=1,2,3$. Then, analogously to the parametric case ([6], Ch. 2), the posterior density of $y=y(\theta)$, with reference to (1) and under $\pi(\cdot)$ as in (2), can be expressed as

$$\begin{aligned}\pi^*(y|X) &= \phi(y)[1+n^{-1/2}(R_1y+R_3y^3)\\ &\quad +n^{-1}\{R_2(y^2-1)+R_4(y^4-3)+R_6(y^6-15)\}]+o_p(n^{-1})\end{aligned} \tag{5}$$

where $X=(X_1,\ldots,X_n)$ and

$$\begin{aligned}R_1 &= a_1(g_3)+m_2^{1/2}\psi_1,\\ R_2 &= b_2(g_3,g_4)+m_2^{1/2}a_1(g_3)\psi_1+\frac{1}{2}m_2(\psi_{11}+\psi_1^2),\\ R_3 &= a_3(g_3),\\ R_4 &= b_4(g_3,g_4)+m_2^{1/2}a_3(g_3)\psi_1,\\ R_6 &= b_6(g_3,g_4).\end{aligned} \tag{6}$$

The propriety of the posterior is assumed here. Let

$$\begin{aligned}u_1 &= R_1+R_3(z^2+2),\\ u_2 &= 2u_1zR_3-\frac{1}{2}u_1^2z+R_2z+R_4(z^3+3z)+R_6(z^5+5z^3+15z),\end{aligned} \tag{7}$$

where $z$ is the $(1-\alpha)$th quantile of a standard normal variate. As in [7], recalling that $y=(n/m_2)^{1/2}(\theta-\bar{X})$, then it follows from (5) that the $(1-\alpha)$th posterior quantile of $\theta$ can be approximated by

$$\theta_1^{(1-\alpha)}(\pi,X)=\bar{X}+(m_2/n)^{1/2}(z+n^{-1/2}u_1), \tag{8}$$

or

$$\theta_2^{(1-\alpha)}(\pi,X)=\bar{X}+(m_2/n)^{1/2}(z+n^{-1/2}u_1+n^{-1}u_2), \tag{9}$$

with posterior coverage error $o_p(n^{-1/2})$ or $o_p(n^{-1})$ respectively.

## 3. Frequentist coverage

### 3.1. Calculation of frequentist coverage

We next study the frequentist coverage of the interval $(-\infty,\theta_2^{(1-\alpha)}(\pi,X)]$. The steps are similar to those in [7] but more involved because of possible data-dependence of the prior; for instance, additional terms appear in the expressions for $W_1$ and $k_2$ in (15) and (21) below.

With $P$ representing the frequentist probability, by (9) and the definition of $y$, the frequentist coverage is given by

$$P\{\theta\le\theta_2^{(1-\alpha)}(\pi,X)\}=P(y\le z+n^{-1/2}u_1+n^{-1}u_2). \tag{10}$$

In order to obtain an expression for the above with margin of error $o(n^{-1})$, we need stochastic expansions for $y$, $u_1$ and $u_2$. To this end, let $E$ denote expectation for fixed $\theta$ and write $\sigma^2 = E(X_i - \theta)^2$, $Z_i = (X_i - \theta)/\sigma$, $\beta_s = E(Z_i^s)$, $1 \le s \le 4$, and $A_s = n^{-1/2}\sum_{i=1}^n (Z_i^s - \beta_s)$, $s = 1, 2, 3$. Then, as in [7],

$$y = -A_1 + \frac{1}{2}n^{-1/2}A_1A_2 - n^{-1}(\frac{1}{2}A_1^3 + \frac{3}{8}A_1A_2^2) + o_p(n^{-1}). \tag{11}$$

Turning next to $u_1$ and $u_2$, we note from (7) that the randomness of these quantities is only due to the $R_i$. For each $i$, let $R_{i0}$ be obtained from $R_i$ in (6) replacing $\bar{X}$, $m_2$, $g_3$ and $g_4$ therein by the corresponding population parameters $\theta$, $\sigma^2$, $\beta_3$ and $\beta_4$ respectively, i.e.,

$$\begin{aligned} R_{10} &= a_1(\beta_3) + \sigma\psi_1^{(0)}, \\ R_{20} &= b_2(\beta_3, \beta_4) + \sigma a_1(\beta_3)\psi_1^{(0)} + \frac{1}{2}\sigma^2[\psi_{11}^{(0)} + \{\psi_1^{(0)}\}^2], \\ R_{30} &= a_3(\beta_3), \\ R_{40} &= b_4(\beta_3, \beta_4) + \sigma a_3(\beta_3)\psi_1^{(0)}, \\ R_{60} &= b_6(\beta_3, \beta_4), \end{aligned} \tag{12}$$

where $\psi_i^{(0)} = \psi_i(\theta, \sigma^2, \beta_3)$, $\psi_{ij}^{(0)} = \psi_{ij}(\theta, \sigma^2, \beta_3)$, $i, j = 1, 2, 3$. Since

$$\bar{X} = \theta + n^{-1/2}\sigma A_1, \quad m_2 = \sigma^2(1 + n^{-1/2}A_2) + o_p(n^{-1/2}), \tag{13}$$

and

$$g_3 = \beta_3 + n^{-1/2}(A_3 - 3A_1 - \frac{3}{2}\beta_3A_2) + o_p(n^{-1/2}), \quad g_4 = \beta_4 + o_p(1), \tag{14}$$

from (6) we get $R_i = R_{i0} + n^{-1/2}W_i + o_p(n^{-1/2})$, $i = 1, 3$, and $R_i = R_{i0} + o_p(1)$, $i = 2, 4, 6$, where

$$\begin{aligned} W_1 &= \sigma^2\psi_{11}^{(0)}A_1 + \sigma\{\frac{1}{2}\psi_1^{(0)} + \sigma^2\psi_{12}^{(0)}\}A_2 \\ &\quad + \{a_1'(\beta_3) + \sigma\psi_{13}^{(0)}\}(A_3 - 3A_1 - \frac{3}{2}\beta_3A_2), \end{aligned} \tag{15}$$

$$W_3 = a_3'(\beta_3)(A_3 - 3A_1 - \frac{3}{2}\beta_3A_2), \tag{16}$$

and $a_i'(\cdot)$ is the derivative of $a_i(\cdot)$. From (7), it is now evident that

$$u_1 = u_{10} + n^{-1/2}\{W_1 + W_3(z^2 + 2)\} + o_p(n^{-1/2}), u_2 = u_{20} + o_p(1),$$

where the leading terms

$$u_{10} = R_{10} + R_{30}(z^2 + 2), \tag{17}$$

and

$$u_{20} = 2u_{10}zR_{30} - \frac{1}{2}u_{10}^2z + R_{20}z + R_{40}(z^3 + 3z) + R_{60}(z^5 + 5z^3 + 15z), \tag{18}$$

are simply counterparts of (7) with the $R_i$ there replaced by the $R_{i0}$.

From (10) and the stochastic expansions for $u_1$ and $u_2$ as indicated above,

$$P\{\theta \le \theta_2^{(1-\alpha)}(\pi, X)\} = P(\tilde{y} \le z + n^{-1/2}u_{10} + n^{-1}u_{20}) + o(n^{-1}), \tag{19}$$

where $\tilde{y} = y - n^{-1}\{W_1 + W_3(z^2+2)\}$. By (11), (15) and (16), the first four approximate cumulants of $\tilde{y}$ are given by

$$\begin{aligned} K_{1n} &= n^{-1/2}k_1 + o(n^{-1}), \\ K_{2n} &= 1 + n^{-1}k_2 + o(n^{-1}), \\ K_{3n} &= n^{-1/2}k_3 + o(n^{-1}), \\ K_{4n} &= n^{-1}k_4 + o(n^{-1}), \end{aligned}$$

where

$$k_1 = \frac{1}{2}\beta_3, \qquad k_3 = 2\beta_3, \qquad k_4 = 12 + 12\beta_3^2 - 2\beta_4 \tag{20}$$

and

$$\begin{aligned} k_2 \;=\; & 3 + \frac{7}{4}\beta_3^2 + 2\{a_1'(\beta_3) + a_3'(\beta_3)(z^2+2) + \sigma\psi_{13}^{(0)}\}(\beta_4 - 3 - \frac{3}{2}\beta_3^2) \\ & + 2\sigma^2\psi_{11}^{(0)} + 2\beta_3\sigma\{\frac{1}{2}\psi_1^{(0)} + \sigma^2\psi_{12}^{(0)}\}. \end{aligned} \tag{21}$$

The fact that $W_1$ and $W_3$ are linear in the $A_i$ facilitates the derivation of (20) and (21). From (19), consideration of an Edgeworth expansion for $\tilde{y}$ now yields

$$P\{\theta \le \theta_2^{(1-\alpha)}(\pi, X)\} = 1 - \alpha + (n^{-1/2}\Delta_1 + n^{-1}\Delta_2)\phi(z) + o(n^{-1}), \tag{22}$$

with

$$\Delta_1 = u_{10} - k_1 - \frac{1}{6}k_3(z^2 - 1) \tag{23}$$

and

$$\begin{aligned} \Delta_2 \;=\; & u_{20} - \frac{1}{2}u_{10}^2 z + zu_{10}\{k_1 + \frac{1}{6}k_3(z^2-3)\} - \frac{1}{2}(k_2 + k_1^2)z \\ & -(\frac{1}{24}k_4 + \frac{1}{6}k_1k_3)(z^3 - 3z) - \frac{1}{72}k_3^2(z^5 - 10z^3 + 15z). \end{aligned} \tag{24}$$

### *3.2. Probability matching conditions*

The frequentist coverage in (22) equals $1 - \alpha + o(n^{-1/2})$ if and only if $\Delta_1 = 0$ identically in $z$ and the population parameters. Since by (12), (17), (20) and (23),

$$\Delta_1 = a_1(\beta_3) + 2a_3(\beta_3) + \sigma\psi_1^{(0)} - \frac{1}{6}\beta_3 + \{a_3(\beta_3) - \frac{1}{3}\beta_3\}z^2,$$

recalling the definition of $\psi_1^{(0)}$, it is clear that the above happens if and only if

$$a_3(\beta_3) = \frac{1}{3}\beta_3 \quad \text{and} \quad \psi(\theta, \sigma^2, \beta_3) = h(\sigma^2, \beta_3) - \theta\sigma^{-1}\{a_1(\beta_3) + \frac{1}{2}\beta_3\}, \tag{25}$$

$h(\cdot)$ being any smooth function of $\sigma^2$ and $\beta_3$. Note that the first condition in (25) is on the empirical-type likelihood whereas the second condition concerns the prior. Indeed, with $\psi(\cdot)$ as in (25), it is easily seen from (2) that the specific choice of $h(\cdot)$ has no influence on the posterior. Hence, hereafter, we take $h(\sigma^2, \beta_3) = 0$, i.e.,

$$\pi(\theta) = \exp[-\theta m_2^{-1/2}\{a_1(g_3) + \frac{1}{2}g_3\}], \tag{26}$$

by (2) and (25), and continue with (25) to obtain further conditions that arise when one wishes to work with margin of error $o(n^{-1})$. Observe that the prior in (26) is actually of the form (4) that motivated the class (2) of priors. From (12), (17), (18), (20), (21) and (24), after considerable algebra, it can be seen that under (25), $\Delta_2 = C_1 z + C_3 z^3 + C_5 z^5$, where

$$\begin{aligned} C_1 &= b_2(\beta_3,\beta_4) + 3b_4(\beta_3,\beta_4) + 15b_6(\beta_3,\beta_4) - \frac{1}{2}\{a_1(\beta_3)\}^2 - \beta_3 a_1(\beta_3) \\ &\quad + \frac{1}{24}\beta_3^2 - \frac{5}{12}\beta_4 + \frac{1}{2}, \\ C_3 &= b_4(\beta_3,\beta_4) + 5b_6(\beta_3,\beta_4) - \frac{1}{3}\beta_3 a_1(\beta_3) + \frac{2}{9}\beta_3^2 - \frac{1}{4}\beta_4 + \frac{1}{2}, \\ C_5 &= b_6(\beta_3,\beta_4) - \frac{1}{18}\beta_3^2. \end{aligned}$$

Thus, under (25), the frequentist coverage in (22) equals $1-\alpha+o(n^{-1})$ for every $z$ and every possible $\theta$, $\sigma^2$, $\beta_3$ and $\beta_4$ if and only if $C_1$, $C_3$ and $C_5$ vanish identically in the population parameters, i.e., if and only if

$$\begin{aligned} b_2(\beta_3,\beta_4) &= \frac{1}{2}\{a_1(\beta_3)\}^2 + \frac{5}{8}\beta_3^2 - \frac{1}{3}\beta_4 + 1, \\ b_4(\beta_3,\beta_4) &= \frac{1}{3}\beta_3 a_1(\beta_3) - \frac{1}{2}\beta_3^2 + \frac{1}{4}\beta_4 - \frac{1}{2}, \\ b_6(\beta_3,\beta_4) &= \frac{1}{18}\beta_3^2. \end{aligned} \tag{27}$$

The conditions in (27) are again on the likelihood.

## 4. Implications

We now examine some major subclasses of empirical-type likelihoods in the light of the conditions obtained in the last section. These are (i) likelihoods arising from empirical discrepancy statistics ([5], Section 1) and hence from Cressie–Read discrepancy statistics [1], (ii) generalized empirical likelihoods [13], and (iii) generalized empirical exponential family likelihoods ([5], Section 4). As noted in [7] all these belong to the general class (1). Moreover, it can be shown that the forms of the $a_i(\cdot)$ and $b_i(\cdot)$ for the subclasses (i)–(iii) are as in Table 1. In this table, $\tau_3$, $\tau_4$, $\gamma_3$, $\gamma_4$ and $\mu$ are constants that depend on the particular likelihood. The usual empirical likelihood belongs to each of (i)–(iii) with

$$\begin{gathered} a_1(g_3) = 0, \qquad a_3(g_3) = \frac{1}{3}g_3 \\ b_0(g_3,g_4) = b_2(g_3,g_4) = 0, \\ b_4(g_3,g_4) = \frac{1}{4}g_4 - \frac{1}{2}(g_3^2+1), \qquad b_6(g_3,g_4) = \frac{1}{18}g_3^2, \end{gathered} \tag{28}$$

while Schennach's Bayesian exponentially tilted empirical likelihood [18] belongs to (iii) with $\mu = \frac{1}{8}$. From (25) and Table 1, it is clear that each likelihood in the subclass (iii) admits, with margin of error $o(n^{-1/2})$, a data-dependent probability matching prior of the form (2) for posterior quantiles. Since each of these likelihoods has $a_1(g_3) = 0$, by (26), such a prior is given by

$$\pi(\theta) = \exp(-\frac{1}{2}\theta m_2^{-1/2} g_3). \tag{29}$$

TABLE 1
*Forms of the $a_i(\cdot)$ and $b_i(\cdot)$ for the subclasses* (i)–(iii)

| Subclass | $a_1(g_3)$ | $a_3(g_3)$ | $b_0(g_3,g_4)$ | $b_2(g_3,g_4)$ | $b_4(g_3,g_4)$ | $b_6(g_3,g_4)$ |
|---|---|---|---|---|---|---|
| (i) | 0 | $\tau_3 g_3$ | 0 | 0 | $\tau_4 g_4 - \frac{9}{2}\tau_3^2(g_3^2+1)$ | $\frac{1}{2}\tau_3^2 g_3^2$ |
| (ii) | 0 | $\gamma_3 g_3$ | 0 | 0 | $\gamma_4 g_4 - \frac{9}{2}\gamma_3^2 g_3^2 - 3\gamma_3 + \frac{1}{2}$ | $\frac{1}{2}\gamma_3^2 g_3^2$ |
| (iii) | 0 | $\frac{1}{3}g_3$ | 0 | 0 | $\mu g_4 - (\mu + \frac{1}{4})(g_3^2+1)$ | $\frac{1}{18}g_3^2$ |

TABLE 2
*Simulation results on the frequentist coverage of $(-\infty, \theta_1^{(1-\alpha)}(\pi, X)]$ for generalized empirical exponential family likelihoods, with $\pi(\theta)$ as in* (29)

| | | Sample size | | | | | Sample size | | | |
|---|---|---|---|---|---|---|---|---|---|---|
| **Distribution** | $1-\alpha$ | **8** | **12** | **16** | **20** | $1-\alpha$ | **8** | **12** | **16** | **20** |
| Normal(0,1) | 0.95 | 0.912 | 0.928 | 0.933 | 0.938 | 0.10 | 0.138 | 0.123 | 0.119 | 0.114 |
| | 0.90 | 0.863 | 0.877 | 0.884 | 0.886 | 0.05 | 0.088 | 0.074 | 0.069 | 0.064 |
| Uniform(0,1) | 0.95 | 0.934 | 0.944 | 0.946 | 0.949 | 0.10 | 0.112 | 0.106 | 0.102 | 0.102 |
| | 0.90 | 0.887 | 0.896 | 0.897 | 0.898 | 0.05 | 0.067 | 0.056 | 0.051 | 0.051 |
| Beta(1,2) | 0.95 | 0.910 | 0.928 | 0.936 | 0.938 | 0.10 | 0.110 | 0.108 | 0.106 | 0.104 |
| | 0.90 | 0.861 | 0.880 | 0.888 | 0.890 | 0.05 | 0.061 | 0.055 | 0.055 | 0.054 |
| Exponential(1) | 0.95 | 0.850 | 0.878 | 0.898 | 0.906 | 0.10 | 0.111 | 0.113 | 0.111 | 0.111 |
| | 0.90 | 0.798 | 0.827 | 0.845 | 0.854 | 0.05 | 0.063 | 0.064 | 0.063 | 0.061 |
| Rayleigh(1) | 0.95 | 0.900 | 0.918 | 0.928 | 0.931 | 0.10 | 0.119 | 0.113 | 0.111 | 0.108 |
| | 0.90 | 0.849 | 0.868 | 0.878 | 0.880 | 0.05 | 0.070 | 0.064 | 0.060 | 0.056 |

From Table 1, it is also clear that none of the likelihoods in the subclasses (i)–(iii) meets the first condition in (27) because $a_1(\cdot) = b_2(\cdot) = 0$ for any such likelihood. Thus, even with possibly data-dependent priors of the form (2), none of them allows frequentist validity of posterior quantiles with margin of error $o(n^{-1})$. In Section 5, it will be seen that at least for the usual empirical likelihood this difficulty can be resolved by considering more elaborate data-dependent priors.

Before addressing this issue, we present some simulation results to indicate the finite sample implications of the aforesaid probability matching property of the prior (29) for the subclass (iii) of generalized empirical exponential family likelihoods. Since this matching holds with margin of error $o(n^{-1/2})$, it makes sense to study the simulated coverage of the interval $(-\infty, \theta_1^{(1-\alpha)}(\pi, X)]$ in this context, where $\theta_1^{(1-\alpha)}(\pi, X)$ approximates the $(1-\alpha)$th posterior quantile of $\theta$ with coverage error $o_p(n^{-1/2})$; see (8). For any likelihood in (iii), it can be seen from (6), (7), (8) and Table 1 that $\theta_1^{(1-\alpha)}(\pi, X) = \bar{X} + (m_2/n)^{1/2}\{z + \frac{1}{6}n^{-1/2}g_3(2z^2+1)\}$, under (29). The simulation results, each based on 10000 simulations are presented in Table 2. Five distributions for the population along with four choices of $1-\alpha$, namely $1-\alpha$=0.95, 0.90, 0.10 and 0.05, are considered. In all cases, except for the exponential distribution in the right tail, the convergence to the desired frequentist coverage turns out to be reasonably fast. Thus the asymptotic results, even with margin of error $o(n^{-1/2})$, are well-reflected in finite samples.

## 5. More elaborate data-dependent priors

Observe that the prior in (26) is equivalent to

$$\pi(\theta) = \exp[-(\theta - \bar{X})m_2^{-1/2}\{a_1(g_3) + \frac{1}{2}g_3\}], \tag{30}$$

in the sense that they both lead to the same posterior; cf. (3) and (4). This motivates us to consider possibly data-dependent priors of the form

$$\pi(\theta) = \exp\{(\theta - \bar{X}) m_2^{-1/2} \chi(g_3) + \frac{1}{2}(\theta - \bar{X})^2 m_2^{-1} \lambda(g_3, g_4)\}, \tag{31}$$

where $\chi(\cdot)$ and $\lambda(\cdot)$ are smooth functions. Since (31) can possibly involve $\bar{X}$ and $g_4$ in addition to $m_2$ and $g_3$, it is more elaborate than (2). The introduction of $\lambda(\cdot)$ in (31) aims at taking care of the population kurtosis.

With reference to any empirical-type likelihood in the general class (1), if one considers the posterior quantiles of $\theta$ under (31), then algebra similar to but heavier than that in Sections 2 and 3 reveals the following:

(a) Frequentist validity of the posterior quantiles holds with margin of error $o(n^{-1/2})$ if and only if

$$a_3(\beta_3) = \frac{1}{3}\beta_3 \text{ and } \chi(\beta_3) = -\{a_1(\beta_3) + \frac{1}{2}\beta_3\}. \tag{32}$$

(b) Frequentist validity of the posterior quantiles holds with margin of error $o(n^{-1})$ if and only if in addition

$$\begin{aligned} b_4(\beta_3, \beta_4) &= \frac{1}{3}\beta_3 a_1(\beta_3) - \frac{1}{2}\beta_3^2 + \frac{1}{4}\beta_4 - \frac{1}{2}, \\ b_6(\beta_3, \beta_4) &= \frac{1}{18}\beta_3^2, \\ \lambda(\beta_3, \beta_4) &= \{a_1(\beta_3)\}^2 - 2b_2(\beta_3, \beta_4) + \frac{5}{4}\beta_3^2 - \frac{2}{3}\beta_4 + 2. \end{aligned} \tag{33}$$

A comparison between (25), (27) and (32), (33) shows that the last condition in (33) is new and this helps. The first condition in (32) as well as the first two conditions in (33) are on the empirical-type likelihood. From Table 1, it can be seen that any likelihood in the subclasses (i)–(iii) meets these three conditions if and only if the associated $a_i(\cdot)$ and $b_i(\cdot)$ are given by (28), which corresponds to the usual empirical likelihood. Furthermore, if (28) holds then the last conditions in (32) and (33) yield $\chi(\beta_3) = -\frac{1}{2}\beta_3$ and $\lambda(\beta_3, \beta_4) = \frac{5}{4}\beta_3^2 - \frac{2}{3}\beta_4 + 2$, so that by (31), the data-dependent prior

$$\pi(\theta) = \exp\{-\frac{1}{2}(\theta - \bar{X}) m_2^{-1/2} g_3 + \frac{1}{2}(\theta - \bar{X})^2 m_2^{-1} (\frac{5}{4} g_3^2 - \frac{2}{3} g_4 + 2)\}, \tag{34}$$

ensures frequentist validity of the posterior quantiles with margin of error $o(n^{-1})$.

A comparison between the prior just obtained and the one shown in (29) is in order. The one in (29) leads to probability matching to a lower order of accuracy but at the same time enjoys the merit of being much simpler. Moreover, as the simulation results reveal, it performs quite well in finite samples. Thus a choice between the two is essentially a matter of taste. If one wishes to work with a simple prior then the one in (29) is recommended. On the other hand, if a premium is put on higher order accuracy from asymptotic considerations, then the one obtained in this section appears to be more attractive.

The connection with a result in [7] is worth noting at this stage. They worked with data-free priors and showed that a likelihood in the general class (1) admits a probability matching prior, with margin of error $o(n^{-1})$, for posterior quantiles if

and only if

$$a_1(\beta_3) = -\frac{1}{2}\beta_3, \qquad a_3(\beta_3) = \frac{1}{3}\beta_3, \qquad b_2(\beta_3, \beta_4) = \frac{3}{4}\beta_3^2 - \frac{1}{3}\beta_4 + 1,$$
$$b_4(\beta_3, \beta_4) = \frac{1}{4}\beta_4 - \frac{2}{3}\beta_3^2 - \frac{1}{2}, \qquad b_6(\beta_3, \beta_4) = \frac{1}{18}\beta_3^2.$$

With the $a_i(\cdot)$ and $b_i(\cdot)$ as above, the conditions in (32) and (33) are met if and only if $\chi(\beta_3) = \lambda(\beta_3, \beta_4) = 0$. In conjunction with (31), one thus gets the flat prior and this agrees with the findings in [7].

## 6. Concluding remarks

The results in Subsection 3.2 show that if $a_3(\beta_3) = \frac{1}{3}\beta_3$ then the prior in (26) ensures frequentist validity of posterior quantiles with margin of error $o(n^{-1/2})$. It is satisfying to note that under the present assumption on the existence of the first four population moments, this margin of error is actually $O(n^{-1})$; vide (22). Similarly, if the existence of the first five population moments is assumed, then under the prior in (34), the frequentist validity of posterior quantiles arising from the usual empirical likelihood holds actually with margin of error $O(n^{-3/2})$.

By (13) and (14), the data-dependent probability matching prior in (34) satisfies $\pi(\theta) = 1 + o_p(1)$ . The same holds for the prior in (30) which is equivalent to the probability matching prior in (26) in the sense of yielding the same posterior. Thus these data-dependent priors approximate the flat prior which is a natural data-free prior in the present context; cf. [7] and a well-known result for fully parametric location models ([6], Chapter 2).

There is scope for extending the present results in several directions. For example, if instead of posterior quantiles, interest lies in the highest posterior density regions, then the findings in [7] show that none of the standard empirical-type likelihoods proposed in the literature, including the usual empirical likelihood, admits a data-free probability matching prior in a higher order asymptotic sense. A natural question is whether consideration of data-dependent priors, such as those of the form (31), can yield more positive results. Another important issue concerns the role of data-dependent priors in the multivariate case with vector $\theta$. Then the algebra will be rather complicated because the empirical-type likelihoods as well as the data-dependent priors will involve multivariate pure and mixed moments in place of $m_2$, $g_3$ and $g_4$, and one would need to consider multivariate Edgeworth expansions for the frequentist calculations. It is hoped that the present work will generate interest in these directions.

**Acknowledgments.** Thanks are due to a referee and the editors for their very constructive suggestions. This work was supported by a grant from the Center for Management and Development Studies, Indian Institute of Management Calcutta.